\documentclass[11pt]{elsarticle}

\usepackage{amssymb}
\usepackage{xcolor}
\usepackage{mathrsfs}
\usepackage{amscd}
\usepackage{appendix}
\usepackage{geometry}
\usepackage{setspace}
\usepackage{amsfonts}
\usepackage{amsmath}
\usepackage{lineno,hyperref}
\usepackage{tikz}
\usetikzlibrary{matrix,arrows}
\newtheorem{definition}{Definition}[section]
\newtheorem{theorem}[definition]{Theorem} 
\newtheorem{lemma}[definition]{Lemma}     

\newtheorem{definition 1}{Definition 1}
\newtheorem{definition 2}{Definition 2}
\newtheorem{remark}[definition]{Remark}

\numberwithin{equation}{section}
\def\Proof{\\\noindent \it Proof.\ \ \rm}
\def\qedbox{$\rlap{$\sqcap$}\sqcup$}
\biboptions{numbers,sort&compress} 
\def\Proof{
\noindent \it Proof.\ \ \rm}
\def\qedbox{$\rlap{$\sqcap$}\sqcup$}

\newcommand{\Z}{{\mathbb Z}}

\journal{arXiv}

\begin{document}

\begin{frontmatter}

\title{Anderson localization for CMV matrices with Verblunsky coefficients defined by the hyperbolic toral automorphism}

\author[]{Yanxue Lin}
\ead{yanxuelin@aliyun.com}

\author[]{Shuzheng Guo\corref{cor1}}
\ead{guoshuzheng@ouc.edu.cn}

\author[]{Daxiong Piao \corref{cor1}}
\ead{dxpiao@ouc.edu.cn}

\address{School of Mathematical Sciences,  Ocean University of China, Qingdao 266100, P.R.China}

\cortext[cor1]{Corresponding authors}

\begin{abstract}
In this paper, we prove the large deviation estimates and Anderson localization for CMV matrices on $\ell^2(\mathbb{Z}_+)$ with Verblunsky coefficients defined dynamically by the hyperbolic toral automorphism. Part of positivity results on the Lyapunov exponents of Chulaevsky-Spencer \cite{C-S95} and Anderson localization results of Bourgain-Schlag \cite{Bourgain02} on Schr\"{o}dinger operators with strongly mixing potentials are extended to CMV matrices.
\end{abstract}
\medskip

\end{frontmatter}



\section{Introduction}

Due to the significance both in physics and mathematics, Anderson localization(AL) problems for Schr\"{o}dinger operators has been widely studied and fruitful results are obtained, see for examples \cite{Bour07, Bourgain03, Bourgain01, BGS02, Bourgain02, ChuSinai89,Jito99, Klein05} and references therein. It is well known that orthogonal polynomials and 1-dimensional discrete Schr\"{o}dinger operator are intimately related (see \cite{B.Simon}).

 In 2018, Wang and Damanik \cite{Wang01} proved AL for the quasi-periodic CMV matrices with Verblunsky coefficients defined by the shift for almost all frequencies, which established the CMV matrix analog of a result proved by Bourgain and Goldstein \cite{Bourgain01}. By following the method of \cite{Bourgain01}, Cedzich and Werner \cite{CW21}  proved AL for one-dimensional quantum walks placed into homogenous electric fields, which implies AL for CMV matrices with a particular choice of skew-shift Verblunsky coefficients for almost all frequencies. Recently, motivated by \cite{Wang01} and \cite{Bourgain02}, we \cite{LPG} considered  AL problems for the CMV matrices with Verblunsky coefficients defined by the skew-shift, and obtained positivity of the Lyapunov exponent and AL for this kind of CMV matrices, which is the CMV matrix analog of the results proved by Bourgain, Goldstein and Schlag for Schr\"{o}dinger operators \cite{Bourgain03}. In 2007, Nguyen \cite{T07} proved positivity of the Lyapunov exponent for CMV matrix defined by expanding map of the unit circle, the large deviation theorem and H\"{o}lder continuity of the Lyapunov exponent, which is a partial CMV generalization of the first main result of \cite{Bourgain02}.

In this paper, we are concerned with the positivity of the Lyapunov exponents and AL for the CMV matrices with Verblunsky coefficients generated by the hyperbolic toral automorphism. Hyperbolic toral automorphism is a strong mixing map on 2-torus; see for example, the paper \cite{C-S95} and the book \cite{Brin}. We shall show that a result of Chulaevski-Spencer \cite{C-S95} concerning the positivity of the Lyapunov exponent in the Schr\"{o}dinger case extends to the OPUC case; see Theorem \ref{mainth02}. And following the formalism of  Bourgain-Schlag \cite{Bourgain02} for the Schr\"{o}dinger operator with potential defined by strong mixing, we show that the corresponding CMV matrix defined on the half-line $\ell^{2}(\mathbb{Z}_+)$ has pure point spectrum in some intervals of the unit circle $\partial \mathbb{D} = \{ z \in \mathbb{C} : |z| = 1 \}$, and the corresponding eigenfunctions is  exponentially decaying, i.e. Anderson localization; see Theorem \ref{mainth01}.


The rest of this paper is structured as follows. We describe the setting and state the main results in Section 2. The large deviation estimate and positivity of Lyapunov exponents are shown in Section 3. The corresponding Green's function estimates are in Section 4 and the main Anderson localization result is then proved in Section 5.

\section{Setting and Results }
In this section, we describe the setting in which we work and state the main results. One may consult \cite{B.Simon,B.Simon-2} for the general background.

Consider the \emph{CMV matrix}, which is a pentadiagonal unitary operator on $\ell^2(\Z_+)$ with a repeating $2\times 4$ block structure of the form
\begin{equation}{\label{1.2}}
\mathcal{C}=\left(
\begin{matrix}
\overline{\alpha}_0&\overline{\alpha}_1\rho_{0}&\rho_1\rho_0&0&0&\cdots&\\
\rho_0&-\overline{\alpha}_1\alpha_{0}&-\rho_1\alpha_0&0&0&\cdots&\\
0&\overline{\alpha}_2\rho_{1}&-\overline{\alpha}_2\alpha_{1}&\overline{\alpha}_3\rho_2&\rho_3\rho_2&\cdots&\\
0&\rho_2\rho_{1}&-\rho_2\alpha_{1}&-\overline{\alpha}_3\alpha_2&-\rho_3\alpha_2&\cdots&\\
0&0&0&\overline{\alpha}_4\rho_3&-\overline{\alpha}_4\alpha_3&\cdots&\\
\cdots&\cdots&\cdots&\cdots&\cdots&\cdots&
\end{matrix}
\right),
\end{equation}
where $\alpha$'s are the \emph{Verblunsky coefficients} which  belong to $\mathbb{D}=\{z\in\mathbb{C}:|z|<1\}$, $\rho_n=\sqrt{1-|\alpha_n|^2}$.

\medskip

Suppose that $\mu$ is a non-trivial (i.e., not finitely supported) probability measure on the unit circle $\partial \mathbb{D} = \{ z \in \mathbb{C} : |z| = 1 \}$. By the non-triviality assumption, the functions $1, z, z^2, \cdots$ are linearly independent in the Hilbert space $\mathcal{H} = L^2(\partial\mathbb{D}, d\mu)$, and hence one can form, by the Gram-Schmidt procedure, the \emph{monic orthogonal polynomials} $\Phi_n(z)$, whose Szeg\H{o} dual is defined by $\Phi_n^{*} = z^n\overline{\Phi_n({1}/{\overline{z}})}$.
The monic orthogonal polynomials obey the  \emph{Szeg\H{o} recurrence}, given by
\begin{equation}\label{eq01}
\Phi_{n+1}(z) = z \Phi_n(z) - \overline{\alpha}_n \Phi_n^*(z), \qquad \textrm{ for } n\in \Z_+.
\end{equation}

\medskip

Consider instead the orthonormal polynomials
$$\varphi_{n}(z)=\frac{\Phi_n(z)}{\|\Phi_n(z)\|_{\mu}},$$
where $\|\cdot\|_{\mu}$ is the norm of $L^{2}(\partial\mathbb{D},d\mu)$. The Szeg\H{o} recurrence  \eqref{eq01} is equivalent to the following one,
\begin{equation}\label{eq01b}
\rho_n \varphi_{n+1}(z)  = z \varphi_{n}(z)  - \overline{\alpha}_n \varphi^*_{n}(z), \textrm{ for } n\in \Z_+.
\end{equation}


The Szeg\H{o} recurrence can be written in a matrix form as follows:
\begin{equation}\label{e.Szego01}
\left(
\begin{matrix}
\varphi_{n+1}\\
\varphi^{*}_{n+1}
\end{matrix}
\right)
=
S_n(x,y;z)
\left(
\begin{matrix}
\varphi_{n}\\
\varphi^{*}_{n}
\end{matrix}
\right), \textrm{ for } n\in\Z_+.
\end{equation}
where
$S(x,y;z)=
\frac{1}{\rho(x,y)}
\left(
\begin{matrix}
z & -\overline{\alpha}(x,y)\\
-\alpha(x,y)z & 1
\end{matrix}
\right)$.
Since det $S(x,y;z)=z$, we study the determinant 1 matrix instead
\begin{equation*}
A_n(x,y;z)
=
\frac{1}{\sqrt{1-|\alpha_n|^2}}
\left(
\begin{matrix}
\sqrt{z} &\frac{ -\overline{\alpha}_n}{\sqrt{z}}\\
-\alpha_n\sqrt{z} & \frac{1}{\sqrt{z}}
\end{matrix}
\right).
\end{equation*}
which is called the \emph{Szeg\H{o} cocycle map}. Then the $n$-step transfer matrix is defined by
\begin{equation}\label{e.n-tr01}
M_n (x,y;z)= \prod^{0}_{j=n-1}A_{j}(x,y;z),
\end{equation}
which is equal to
\begin{equation}\label{e.transfer01}
M_n(x,y;z)=\frac{1}{2}
\left(
\begin{matrix}
\varphi_{n}(z)+\psi_{n}(z) & \varphi_{n}(z)-\psi_{n}(z)\\
\varphi_{n}^\ast(z)-\psi_{n}^\ast(z)&\varphi_{n}^\ast(z)+\psi_{n}^\ast(z)
\end{matrix}
\right).
\end{equation}
where $\{\psi_n(z)\}_{n=0}^{\infty}$ are the second kind polynomials associated to the $\{\varphi_n(z)\}_{n=0}^{\infty}$, i.e., the set of orthonormal polynomials whose Verblunsky coefficients are given by $\{-\alpha_n\}_{n=0}^{\infty}$.

We have the following expression for the Lyapunov exponent in terms of the $\varphi_{n}(z)$ and $\psi_{n}(z)$ \cite [Proposition 10.5.5]{B.Simon-2}:
\begin{equation}\label{e.Lyapunov3}
L(z)=\lim _{n \rightarrow \infty}\frac{1}{2 n} \log \left(\left|\varphi_{n}(z)\right|^{2}+\left|\psi_{n}(z)\right|^{2}\right).
\end{equation}




\medskip

In this paper, we consider a sequence of Verblunsky coefficients generated by sampling function $\alpha(\cdot,\cdot):\mathbb{T}^2\rightarrow\mathbb{D}$, i.e., $\alpha_{n}(x,y)=\lambda \alpha(A^n(x,y))$, where $(x,y)\in\mathbb{T}^2$, $\mathbb{T}^2=\mathbb{R}^2/2\pi\mathbb{Z}^2$, $\alpha\in C^1(\mathbb{T}^2, \mathbb{D})$ satisfying $\langle\alpha(x,y)\rangle=0$,  $\lambda\in(0,1)$  is a coupling parameter, and $A: \mathbb{T}^2 \rightarrow \mathbb{T}^2$ is a hyperbolic toral automorphism.
We assume that the sampling function $\alpha(x,y)$ satisfies
\begin{equation*}
\int_{\mathbb{T}^{2}}\log(1-|\alpha(x,y)|)dxdy > -\infty.
\end{equation*}
With $\{\alpha_n\}_{n\in\Z_+}$ as above, the CMV matrix $\mathcal{C}$ can be generated dynamically. Before stating the main results, we introduce the following definitions.
\begin{definition}
(Anderson localization) We say $\mathcal{C}$ exhibits \emph{Anderson localization} on $\mathcal{I}\subset\partial\mathbb{D}$, if $\mathcal{C}$ has only pure point spectrum in $\mathcal{I}$ and its eigenfunction $\xi_{n}$ decay exponentially in $n$.
\end{definition}
\begin{definition}
(Generalized eigenvalues and generalized eigenfunctions) We call $z$ a \emph{generalized eigenvalue} of $\mathcal{C}$, if there exists a nonzero, polynomially bounded function $\xi_{n}$ such that $\mathcal{C} \xi=z\xi$, where $\xi=\{\xi_n\}_{n\in\mathbb{Z}_+}$. We call $\xi_{n}$ a \emph{generalized eigenfunction}.
\end{definition}

For $f:\mathbb{T}^2 \rightarrow \mathbb{C}$, define its mean value as $ \langle f\rangle:=\frac{1}{4\pi^2}\int_{\mathbb{T}^2}f(x,y)dxdy$.
Let $F_n=\alpha\circ A^n$. Define the function
\begin{equation}
\mathcal{J}(\eta)=\sum_{n=-\infty}^{\infty}e^{in\eta}\langle \bar{F}_{0}F_{n}\rangle,
\end{equation}
we have the natural definition $\langle\bar{F}_{0}F_{n}\rangle:=\langle\bar{F}_{-n}F_{0}\rangle$ for $n<0$.

Through out this paper, we suppose the mean value of sampling function $\alpha$ is zero, that is $\langle\alpha\rangle=0$.
\begin{remark}

\begin{enumerate}
\item Notice that $A$ is a hyperbolic toral automorphism, $\det(A)=1$, so $\langle F_n\rangle=\langle\alpha\circ A^n\rangle=\langle \alpha\rangle=0$.
\item Due to the conjugate of $e^{i n\eta}\langle\bar{F}_{0}F_{n}\rangle$ and $e^{-in\eta}\langle\bar{F}_{-n}F_{0}\rangle$, it is easy to verify that $\mathcal{J}(\eta)$ is real-valued.
\end{enumerate}
\end{remark}

Now we can state our main results as follows.
\begin{theorem}\label{mainth02}
For $z\in\partial\mathbb{D}\backslash \{1,-1\}$, we have for $\lambda \in (\lambda_0,1)$ with $\lambda_0\in(0,1)$ that
\begin{equation}\label{lyapunov01}
L(z)=\frac{\lambda^2}{2}\mathcal{J}(\eta)+O(\lambda^{3}).
\end{equation}
Moreover, for those $\lambda$ and large $N$, we have the large deviation estimate
\begin{equation}\label{LDT01}
\mathrm{mes}\left[(x,y)\in \mathbb{T}^2:\left| L_N(z)- \frac{\lambda^2}{2}\mathcal{J}(\eta)\right|>\lambda^{3}  \right]<\mathrm{exp}(-C_\lambda N).
\end{equation}
\end{theorem}

\begin{theorem}\label{mainth01}
Fix some small $\delta>0$ and a sufficiently small $\lambda >0$. Let $I_0\subset[e^{i\delta},e^{i(\pi-\delta)}]\bigcup$ $[e^{i(\pi+\delta)},e^{i(2\pi-\delta)}]$ such that \eqref{e.spectral01} holds on $I_0$. Then there exists a set $\Omega\subset\mathbb{T}^2$ with $\mathrm{mes}(\mathbb{T}^2\setminus\Omega)\leq\varepsilon$, such that for every $(x,y)\in\Omega$, the CMV matrix on $\ell^{2} (\mathbb{Z}_+)$ has pure point spectrum in the interval $I_0$ and the corresponding eigenfunctions decay exponentially.
\end{theorem}


\section{A Large Deviation Estimate}

We use the method developed by Bourgain and Schlag \cite{Bourgain02} to prove Anderson localization, which means that a large deviation estimate will play an important role. We will present the transfer matrix we need in the following sections, along with the Pr\"{u}fer variables, which are essentially polar coordinates for the orthonormal polynomials $\varphi_{n}(z)$.

\subsection{Pr\"{u}fer Variables}
In this subsection, we will write down the analogous variables for OPUC when $z \in \partial\mathbb{D}$ and use them as a tool for spectral analysis, see \cite[section 10.12]{B.Simon}.

Let $z=e^{i\eta}\in\partial\mathbb{D}$ $(\eta \in [0,2\pi))$. Define the corresponding Pr\"{u}fer variables $R_{n}(z)$ and $\theta_{n}(z)$ by
\begin{equation}\label{e.prufer01}
\Phi_{n}(z)=R_{n}(z) \mathrm{exp}[i(n\eta+\theta_{n}(z))],
\end{equation}
where $\theta_{n}$ is real and determined by $|\theta_{n+1}-\theta_{n}|<\pi$. Here $R_{n}(z)>0$, that is
\begin{equation*}
R_{n}(z)=|\Phi_{n}(z)|.
\end{equation*}

\medskip

In this paper, $\varphi_n(z)$ enters more often than $\Phi_{n}(z)$, so we use the following expression equivalent to \eqref{e.prufer01},
\begin{equation*}
\varphi_{n}(z)=r_{n}(z) \exp [i\left(n \eta+\theta_{n}(z)\right)],
\end{equation*}
where $r_{n}(z)=|\varphi_{n}(z)|$. Alternatively, it drops out of $\varphi^{\ast}_{n}(z)$, i.e.,
\begin{equation*}
\varphi^{\ast}_{n}(z)=r_{n}(z)\exp[-i\theta_{n}(z)].
\end{equation*}

\medskip

For $z=e^{i\eta}\in\partial\mathbb{D}$, define $\zeta_{n}(z)=e^{i\left[(n+1) \eta+2 \theta_{n}(z)\right]}$, we have
\begin{equation}\label{e.prufer2}
\frac{r_{n+1}^{2}}{r_{n}^{2}}=\frac{1+\left|\alpha_{n}\right|^{2}-2 \operatorname{Re}\left(\alpha_{n} \zeta_{n}\right)}{1-\left|\alpha_{n}\right|^{2}}\equiv H_{n},
\end{equation}
\begin{equation*}
\exp [-i\left(\theta_{n+1}-\theta_{n}\right)]=\frac{\rho_{n}^{-1}(1-\alpha_{n} \zeta_n)}{H_{n}^{1 / 2}}
\end{equation*}
and
\begin{align}
\zeta_{n+1} &=\frac{z\zeta_{n}\left(1+\left|\alpha_{n}\right|^{2}-2 \operatorname{Re}\left(\alpha_{n} \zeta_{n}\right)\right)}{\left(1-\alpha_{n} \zeta_{n}\right)^{2}} \label{e.prufer4}\\
&=z \zeta_{n}\left[1+2 i \operatorname{Im}\left(\alpha_{n} \zeta_{n}\right)\right]+O\left(\lambda^{2}\right) \label{e.prufer4a}.
\end{align}

Iterating \eqref{e.prufer4a}, we have for $T\sim\log1/\lambda$ to be determined later that
\begin{align}
\zeta_n & =z^T \zeta_{n-T} \prod_{s=1}^T\left[1+2 i \operatorname{Im}\left(\alpha_{n-s} \zeta_{n-s}\right)\right]+O\left(\lambda^2 \log \lambda\right) \notag \\
& =z^T \zeta_{n-T}\left[1+\sum_{s=1}^T 2 i \operatorname{Im}\left(\alpha_{n-s} \zeta_{n-s}\right)\right]+\sum_{n=2}^{T-1} O\left(T^n \lambda^n\right)+O\left(\lambda^2 \log \lambda\right) \notag\\
& =z^T \zeta_{n-T}\left[1+2 i \operatorname{Im}\left(\zeta_{n-T} \sum_{s=1}^T z^{T-s} \alpha_{n-s}\right)\right]+O\left(\lambda^2 \log ^2 1 / \lambda\right) . \label{e.prufer6}
\end{align}

Recall that $F_n=\alpha\circ A^n$. According to \eqref{e.Lyapunov3}, using the Pr\"{u}fer recurrence, we can write
\begin{align}
&\frac{1}{2 N} \log \left|\varphi_{N}(z)\right|^{2}\notag=\frac{1}{2 N} \log r_{N}^{2}(z) \notag\\
&=\frac{1}{2 N} \sum_{n=1}^{N} \log \left[\left(1-\left|\alpha_{n}\right|^{2}\right)^{-1}\right]+\frac{1}{2 N} \sum_{n=1}^{N} \log \left(1+\left|\alpha_{n}\right|^{2}-2 \operatorname{Re}\left(\alpha_{n} \zeta_{n}\right)\right)\notag\\
&=\frac{\lambda^{2}}{2 N} \sum_{n=1}^{N}\left|F_{n}\right|^{2}+\frac{1}{2 N} \sum_{n=1}^{N}\left(\left|\alpha_{n}\right|^{2}-2 \operatorname{Re}\left(\alpha_{n} \zeta_{n}\right)-2 \operatorname{Re}\left(\alpha_{n} \zeta_{n}\right)^{2}\right)+O\left(\lambda^{3}\right)\notag\\
& = \frac{\lambda^{2}}{2 N} \sum_{n=1}^{N}\left|F_{n}\right|^{2} -\frac{\lambda}{N} \sum_{n=1}^{N} \operatorname{Re}\left(\zeta_{n} F_{n}\right)-\frac{\lambda^{2}}{4 N} \sum_{n=1}^{N}\left(\zeta_{n}^{2} F_{n}^{2}+\bar{\zeta}_{n}^{2} \bar{F}_{n}^{2}\right)+O\left(\lambda^{3}\right)\notag\\
&=:I_1+I_2+I_3.\label{e.prufer001}
\end{align}
Similarly for $\frac{1}{2N}\log|\psi_N(z)|^2$ with the replacement $F_n \longmapsto -F_n$ in the above.

\subsection{Large Deviation Estimate}



For any $f \in C^1(\mathbb{R})$ and $h>0$, let
$$
m_f(h)=\sup _{I:|I|=h}\left|\int_I\hspace{-1.05em}- f d y\right|,
$$
where the supremum is taken over all intervals $I \subset \mathbb{T}$ and $\int_I\hspace{-1.05em}- :=\frac{1}{|I|} \int_I$.

The following lemma is contained in \cite{Bourgain02}. For convenience, we only state it here without proof.
\begin{lemma}(\cite[Lemma 10.1]{Bourgain02})\label{le.10.1}
Let $f$, $\left\{g_j\right\}_{j=1}^{\infty} \in C^1(\mathbb{R})$,  $\|f\|_{\infty} \leq 1$, $\sup _j\left\|g_j\right\|_{\infty} \leq 1$. Suppose $\left\|f^{\prime}\right\|_{\infty}+$ $\mathrm{sup} _j\left\|g_j^{\prime}\right\|_{\infty} \leq \widetilde{K}$. Then for any $L \geq 2$,
\begin{equation}\label{line01}
 \sup _{|I|=1} \int_I \exp \left(t \sum_{j=1}^N g_j\left(L^{j-1} y\right) f\left(L^j y\right)\right) d y \leq e^{C t^2 N},
\end{equation}
if
\begin{equation}\label{line02}
 \widetilde{K} L^{-\frac{1}{2}}+m_f(\sqrt{L}) \leq t.
\end{equation}
\end{lemma}

\begin{lemma}\label{le.LDT02}
Fix hyperbolic toral automorphism $A$ on $\mathbb{T}^2$. Let $F\in C^1(\mathbb{T}^2,\mathbb{D})$, $\langle F(x,y)\rangle=0$, and such that $|F(x,y)|\leq1$, $|\nabla\mathrm{Re}(F)|+|\nabla\mathrm{Im}(F)|\leq K$. Then for any $0<\delta<1$ and all $N \geq N_0(A, K, \delta)$,
\begin{equation}\label{e.LDT01}
\mathrm{mes}\left[(x,y)\in\mathbb{T}^2:\left|\frac{1}{N}\sum_{n=1}^{N}F_n(x,y)\right|>\delta\right]<C\mathrm{exp}(-c\delta^2\frac{N}{\log (K/\delta)}),
\end{equation}
where $K>1$, $c$, $C$ are constants depending only on $A$.
\end{lemma}
\Proof
The proof is similar to \cite[Corollary 10.2]{Bourgain02} and \cite[Lemma 3]{T07}, and will be omitted here.
\qedbox

\medskip

Now we estimate $I_1(=\frac{\lambda^{2}}{2 N} \sum_{n=1}^{N}\left|F_{n}\right|^{2})$ in \eqref{e.prufer001}.
Fix some small $\delta>0$, by \eqref{e.LDT01} with $K\sim 1$,
\begin{equation}\label{e.LDT02}
\mathrm{mes}\left[(x,y)\in\mathbb{T}^2:\left|\frac{1}{N}\sum_{n=1}^{N}|F_n(x,y)|^2-\left\langle |F_0|^2\right\rangle\right|>\delta\right]<C\mathrm{exp}(-c\delta^2\frac{N}{|\log \delta|}).
\end{equation}

\medskip

To estimate $I_3(= -\frac{\lambda^{2}}{4 N} \sum_{n=1}^{N}\left(\zeta_{n}^{2} F_{n}^{2}+\overline{\zeta}_{n}^{2} \overline{F}_{n}^{2}\right))$ in \eqref{e.prufer001}, one need to introduce some notations. Let $A\underline{\upsilon}^{+}=\varrho\underline{\upsilon}^{+}$, $A\underline{\upsilon}^{-}=\varrho_{-}\underline{\upsilon}^{-}$, where $|\varrho|>1$, $\varrho_{-}=\pm\frac{1}{\varrho}$ and $\underline{\upsilon}^{+}$, $\underline{\upsilon}^{-}\in\mathbb{T}^2$ are unit vectors, $||\underline{\upsilon}^{+}||=||\underline{\upsilon}^{-}||=1$. From \eqref{e.prufer4}, we have $\zeta_n=V(\zeta_{n-1},F_{n-1})$, where
\begin{equation*}
V(\zeta,F)=\frac{z \zeta\left(1+\lambda^2\left|F\right|^{2}-2\lambda \operatorname{Re}\left(F \zeta\right)\right)}{\left(1-\lambda F \zeta\right)^{2}}.
\end{equation*}
Notice that $|V_\zeta|\leq1+C\lambda$, $|V_F|\leq C\lambda$ for small $\lambda$, and $|F|\leq1$. Thus
\begin{equation}\label{e.dif01}
|\nabla\zeta_n|\leq(1+C\lambda)|\nabla \zeta_{n-1}|+C\lambda\varrho^{n-1}\leq C\lambda\varrho^{n-1}\sum_{\ell=0}^\infty\left(\frac{1+C\lambda}{\varrho}\right)^\ell\leq C\lambda\varrho^{n-1}
\end{equation}
for small enough $\lambda$. Therefore, we have
\begin{equation}\label{e.dif02}
||\nabla\zeta_n||_\infty\leq C\lambda\varrho^{n-1},
\end{equation}
for all $n$.

\begin{lemma}
Let $F\in C^1(\mathbb{T}^2,\mathbb{D})$ satisfying $|\nabla\mathrm{Re}(F)|+|\nabla\mathrm{Im}(F)|\lesssim 1$. For any $0<\lambda<\lambda_0(\alpha,A)$ and $N\geq N_0(\lambda)$,
\begin{equation}\label{e.LDT04}
\mathrm{mes}\left[(x,y)\in\mathbb{T}^2:\left|\frac{1}{N}\sum_{n=1}^{N}\left(F_n^2\zeta_n^2 +\bar{F}_n^2\bar{\zeta}^2_n\right) \right|>\lambda^{3}  \right]\leq e^{-C_\lambda N}.
\end{equation}
\end{lemma}
\Proof
\begin{equation*}
\begin{split}
\left|\frac{1}{N}\sum_{n=1}^{N}\left(F_n^2\zeta_n^2+\bar{F}_n^2\bar{\zeta}^2_n\right)\right|&\leq \left|\frac{1}{N}\sum_{n=1}^{N}F_n^2\zeta_n^2\right|+\left|\frac{1}{N}\sum_{n=1}^{N}\bar{F}_n^2\bar{\zeta}^2_n\right|\\
&=2\left|\frac{1}{N}\sum_{n=1}^{N}F_n^2\zeta_n^2\right|.
\end{split}
\end{equation*}
It suffices to show that
\begin{equation}\label{e.LDT07}
\mathrm{mes}\left[(x,y)\in\mathbb{T}^2:\left|\frac{1}{N}\sum_{n=1}^{N}F_n^2\zeta_n^2\right|>\frac{\lambda^{3}}{2}  \right]\leq e^{-C_\lambda N}.
\end{equation}
Recall that $\zeta_n=e^{i[(n+1)\eta+2\theta_n(z)]}$, where $\eta\in[0,2\pi)$, $z=e^{i\eta}$ and $\theta_n$ is real.
According to \eqref{e.prufer4a}, we have
$$\frac{1}{N}\sum_{n=1}^N\zeta_n=z(\frac{1}{N}\sum_{n=1}^N\zeta_{n})+O(\lambda )+O(\frac{1}{N}),$$
so that
$$(1-z)\left|\frac{1}{N}\sum_{n=1}^N\zeta_n\right|= O(\lambda ),$$
$$(1-z^2)\left|\frac{1}{N}\sum_{n=1}^N\zeta^2_n\right|= O(\lambda ).$$
Therefore, for $z=e^{i\eta}(\neq\pm1)$, we have
$$\left|\frac{1}{N} \sum_{n=1}^N \zeta_n\right|+\left|\frac{1}{N} \sum_{n=1}^N \zeta_n^2\right| =O(\frac{\lambda}{\sin^2\eta})\leq C \lambda.$$

Let $W=F^2-\left\langle F^2\right\rangle$ so that $F_n^2(x,y)-\left\langle F_n^2\right\rangle=W\left(A^n (x,y)\right)=W_n(x,y)$. It suffices to prove \eqref{e.LDT07} with $W_n$ instead of $F_n^2$ up to an error $O(\frac{\lambda^3}{\sin^2\eta})$. Let $\zeta^2_n=\xi_n+i \eta_n$, $W_n(x,y)=P_n+iQ_n$. We have
 $$W_n\zeta^2_n=(P_n \xi_n- Q_n\eta_n)+i(Q_n\xi_n+P_n\eta_n).$$
We will show that
\begin{equation}\label{e.LDT05}
\mathrm{mes}\left[(x,y)\in \mathbb{T}^2: \frac{1}{N}\sum_{n=1}^{N}P_n\xi_n(x,y)>\lambda^3\right]\leq e^{-C_\lambda N}.
\end{equation}
The other cases are similar.

Fix the positive integer $s$. Since $|P_n(x,y)-P_n(\underline{\upsilon}^+\underline{\upsilon}^+\cdot(x,y))|\leq|\nabla P|\varrho^{-n}$, we have
\begin{equation}\label{e.LDT06}
\begin{aligned}
\frac{1}{N} \sum_{n=1}^N P_n(x,y) \xi_n & =\frac{1}{s} \sum_{r=1}^s \frac{1}{N / s} \sum_{j=0}^{[N / s]-1}\left(P_{j s+r} \xi_{j s+r}\right)(x,y)+O\left(\frac{s}{N}\right) \\
& =\frac{1}{s} \sum_{r=1}^s \frac{1}{N / s} \sum_{j=0}^{[N / s]-1} P_{j s+r}\left(\underline{v}^{+} (x,y) \cdot \underline{v}^{+}\right) \xi_{j s+r}(x,y)+O\left(\frac{s}{N}\right) .
\end{aligned}
\end{equation}

Let $f(\widetilde{y})=P(\underline{\upsilon}^+ \widetilde{y})$ and $g_{j,r}(\widetilde{y},\widetilde{y}^{'})=\xi_{js+r}(\varrho^{-s(j-1)-r}\underline{\upsilon}^+ \widetilde{y}+\underline{\upsilon}^-\widetilde{y}^{'})$. Denote $\widetilde{y}=\underline{v}^{+}\cdot(x,y)$, $\widetilde{y}^{'}=\underline{v}^{-}\cdot(x,y)$ as the projection lengths of $(x,y)$ on $\underline{v}^{+}$ and $\underline{v}^{-}$ respectively. By \eqref{e.dif02},
\begin{equation*}
\left|\frac{\partial}{\partial\widetilde{y}} g_{j,r}(\widetilde{y},\widetilde{y}^{'})\right|\leq C\lambda\varrho^s.
\end{equation*}

Let $P_R=\phi_{\frac{1}{R}}\ast P$, where $\phi_{\frac{1}{R}}$ is an $L^1$-normalized bump function with $\mathrm{diam}(\mathrm{supp}(\phi_{\frac{1}{R}}))\thicksim 1/R$. It is obvious that
\begin{align}
||P_R(\underline{\upsilon}^+ \widetilde{y})-P(\underline{\upsilon}^+ \widetilde{y})||_\infty &=\sup_{(x,y)\in\mathbb{T}^2}|\phi_{\frac{1}{R}}\ast P(\underline{\upsilon}^+ \widetilde{y})-P(\underline{\upsilon}^+ \widetilde{y})|\notag\\
&=\sup_{(x,y)\in\mathbb{T}^2}\left|\int_{\mathbb{T}^2}\int_{\mathbb{T}^2}\phi_{\frac{1}{R}}(t_1,t_2)P(\underline{\upsilon}^+ \widetilde{y}-(t_1,t_2))dt_1dt_2-P(\underline{\upsilon}^+ \widetilde{y})\right|\notag\\
&<\sup_{(x,y)\in\mathbb{T}^2}\int\int_{[\mathrm{supp}(\phi_{\frac{1}{R}})]^2}\phi_{\frac{1}{R}}(t_1,t_2)dt_1dt_2 \cdot \sup_{(x,y)\in\mathbb{T}^2}|\nabla P(\underline{\upsilon}^+ \widetilde{y})|\notag\\
&<1\cdot\frac{1}{R^2}.\label{e.LDT+}
\end{align}
Set $\bar{f}(\widetilde{y})=P_R(\underline{v}^{+}\widetilde{y})$. Then $||\bar{f}^{'}||_\infty\lesssim 1$ and for any interval $I\subset \mathbb{R}$,
\begin{equation*}
\begin{split}
hm_{\bar{f}}(h)=\left|\int_I \bar{f}(\widetilde{y})d\widetilde{y}\right|&=\left|\int_I \phi_{\lambda^3}\ast P(\underline{\upsilon}^+ \widetilde{y})d\widetilde{y}\right|\\
&\leq\sum_{|(k_1,k_2)|>0}\left|\widehat{P}(k_1,k_2)\right|\left|\widehat{\phi}_{\lambda^3}(k_1,k_2)\right|\left| \int_I e^{i(k_1,k_2)\cdot\underline{v}^{+}\widetilde{y}}d\widetilde{y}\right|\\
&\leq\sum_{|(k_1,k_2)|>0}\left|\widehat{P}(k_1,k_2)\right|\left|\widehat{\phi}_{\lambda^3}(k_1,k_2)\right| \left\|(k_1,k_2)\cdot\underline{v}^{+}\right\|^{-1}\leq C\lambda^{-12},
\end{split}
\end{equation*}
with $R\thicksim \lambda^{-3}$,
\eqref{e.LDT+} shows that the same result holds for $f(\widetilde{y})$, which is $hm_{f}(h)\leq C\lambda^{-12}$.

According to the exponential Chebyshev's inequality, with large $N$ and  $L=\varrho^s$,
\begin{align}
&\mathrm{mes}\left[(x,y) \in\mathbb{T}^2: \frac{1}{N / s} \sum_{j=0}^{[N / s]-1} P_{j s+r}(x,y) \xi_{j s+r}(x,y)>\lambda^3\right] \notag\\
& \leq \int_{\mathbb{T}^2}\mathrm{exp}\left(t\sum_{j=0}^{[N/s]-1}P_{js+r}\xi_{js+r}-\lambda^3t\frac{N}{s}\right)dxdy\notag\\
& \leq C \exp(-\lambda^3 t N / s) \int_0^1 \sup _{|I|=1} \int_I \exp \left(t \sum_{j=0}^{[N / s]-1} f\left(L^j \widetilde{y}\right) g_{j, r}\left(L^{j-1} \widetilde{y}, \widetilde{y}'\right)\right) d \widetilde{y} d \widetilde{y}^{\prime}\notag \\
& \leq C \exp \left(-\lambda^3 t N / s+C t^2 N / s\right),\label{LDT02}
\end{align}
provided
\begin{equation*}
(1+\lambda\varrho^s)L^{-\frac{1}{2}}+\lambda^{-12}L^{-\frac{1}{2}}\leq t.
\end{equation*}
Here, we used Lemma~\ref{le.10.1} in the third step.
Setting $L=\lambda^{-30}$ and $\varrho^s=\lambda^{-11}$ one has $t\sim\lambda^{3}$. Since \eqref{e.LDT06} and \eqref{LDT02} imply \eqref{e.LDT05} for large $N$, the lemma follows.
\qedbox

\medskip

Next, we estimate $I_2(=-\frac{\lambda}{N}\sum_{n=1}^{N}\mathrm{Re}(\zeta_n F_n))$ in \eqref{e.prufer001}. Fix an integer $T$ so that $\varrho^T \sim\lambda^{-100}$. By using the recurrence \eqref{e.prufer6}, one can get
\begin{align}
I_2= & -\frac{\lambda}{N} \sum_{n=T+1}^N \mathrm{Re}(z^T \zeta_{n-T} F_n )\notag \\
& +\frac{\lambda^2}{N} \sum_{n=T+1}^N \sum_{s=1}^T \mathrm{Re}(z^s \bar{F}_{n-s} F_n) \notag \\
& -\frac{\lambda^2}{N} \sum_{n=T+1}^N \sum_{s=1}^T \mathrm{Re}(z^{2 T-s} \zeta_{n-T}^2 F_{n-s} F_n)+O\left(\lambda^3 \log ^2 1 / \lambda\right)+O(T / N)\notag\\
&=:I_4+I_5+I_6.\label{e.prufer002}
\end{align}

Firstly, we estimate $I_5(=\frac{\lambda^2}{N} \sum_{n=T+1}^N \sum_{s=1}^T \mathrm{Re}(z^s \bar{F}_{n-s} F_n))$ in \eqref{e.prufer002}. By Lemma \ref{le.LDT02} with $K=\varrho^s$, $\delta=\lambda^3$, we have
\begin{equation*}
\mathrm{mes}\left[(x,y)\in\mathbb{T}^2:\left|\frac{1}{N}\sum_{n=1}^{N}\bar{F}_{n-s}F_n-\langle \bar{F}_0F_s\rangle\right|>\lambda^3\right]\leq\mathrm{exp}(-C_\lambda N).
\end{equation*}

Secondly, we estimate $I_4(=-\frac{\lambda}{N} \sum_{n=T+1}^N \mathrm{Re}(z^T \zeta_{n-T} F_n ))$ in \eqref{e.prufer002}.
\begin{lemma}\label{le.LDT03}
For $N>N_0(\lambda, \varrho)$, we have
\begin{equation}\label{e.LDT03}
\mathrm{mes}\left[(x,y)\in\mathbb{T}^2: \left|\frac{1}{N}\sum_{n=T+1}^{N}\mathrm{Re}(F_{n}\zeta_{n-T})\right|>\lambda^3 \right]\leq \mathrm{exp}(-C_\lambda N).
\end{equation}
\end{lemma}

\Proof
Recall that $\zeta_n=\xi_n+i\eta_n$, let $\mathrm{Re}(F_n\zeta_{n-T})=A_n\xi_{n-T} -B_n \eta_{n-T} $. As in \eqref{e.LDT06},
\begin{equation*}
\begin{aligned}
\frac{1}{N} \sum_{n=T+1}^N A_n (x,y)\xi_{n-T}(x,y) & =\frac{1}{T} \sum_{r=1}^T \frac{1}{N / T} \sum_{j=1}^{[N / T]}A_{j T+r}(x,y) \xi_{j T+r}(x,y)+O\left(\frac{T}{N}\right) \\
& =\frac{1}{T} \sum_{r=1}^T \frac{1}{N / T} \sum_{j=1}^{[N / T]} A_{j T+r}\left[\underline{v}^{+} (\underline{v}^{+}\cdot (x,y))\right]\xi_{j T+r}(x,y)+O\left(\frac{T}{N}\right) .
\end{aligned}
\end{equation*}
It suffices to prove that for any $r=1,2,\cdots, T$,
\begin{equation}
\mathrm{mes}\left[ (x,y)\in\mathbb{T}^2:\frac{1}{N / T} \sum_{j=1}^{[N / T]}A_{j T+r} (x,y)\xi_{(j-1) T+r}(x,y)\geq\lambda^3\right]\leq\mathrm{exp}(-C_\lambda N).
\end{equation}
As before, let $f(\widetilde{y})=A(\underline{v}^{+}\widetilde{y})$, $g_{j,r}(\widetilde{y},\widetilde{y}^{'})=\xi_{(j-1)T+r}(\varrho^{-(j-1)T-r}\underline{v}^+\widetilde{y}+\underline{v}^-\widetilde{y}^{'})$. Note that
\begin{equation*}
A_{jT+r}(\varrho^{-r}\underline{v}^{+}\widetilde{y})\xi_{(j-1)T+r}(\varrho^{-r}\underline{v}^+\widetilde{y}+\underline{v}^-\widetilde{y}^{'})=f(L^j\widetilde{y})g_{j,r}(L^{j-1}\widetilde{y},\widetilde{y}^{'})
\end{equation*}
with $L=\varrho^T\thicksim\lambda^{-100}$. Moreover, $|\frac{\partial}{\partial \widetilde{y}}g_{j,r}(\widetilde{y},\widetilde{y}^{'})|\leq C\lambda$, and $hm_f(h)\leq\lambda^{-12}$. Hence,  with $t\sim\lambda^3$,
\begin{equation*}
\begin{aligned}
&\mathrm{mes}\left[(x,y)\in\mathbb{T}^2:\frac{1}{N/T}\sum_{j=1}^{[N/T]}A_{j T+r}(x,y) \xi_{(j-1) T+r}(x,y)>\lambda^3\right]\\
&\leq\int_{\mathbb{T}^2}\mathrm{exp}\left[\left(\sum_{j=1}^{[N/T]}A_{j T+r} \xi_{(j-1) T+r}-\lambda^3 N/T\right)t\right] dxdy\\
&\leq C\mathrm{exp}(-t\lambda^3N/T)\int_{0}^{1}\sup_{|I|=1}\int_I\mathrm{exp}\left(t\sum_{j=1}^{[N/T]}f(L^j\widetilde{y})g_{j,r}(L^{j-1}\widetilde{y},\widetilde{y}^{'}) \right)d\widetilde{y}d\widetilde{y}^{'}\\
&\leq\mathrm{exp}(-t\lambda^3 N/T+Ct^2N/T)\leq\mathrm{exp}(-C_\lambda N),
\end{aligned}
\end{equation*}
where we used Lemma \ref{le.10.1} in the third step. Notice that condition \eqref{line02} is satisfied, since $(1+C\lambda)L^{-\frac{1}{2}}+\lambda^{-12}L^{-\frac{1}{2}}\leq t.$
$\square$

\medskip

Then, we estimate $I_6(=-\frac{\lambda^2}{N} \sum_{n=T+1}^N \sum_{s=1}^T \mathrm{Re}(z^{2 T-s} \zeta_{n-T}^2 F_{n-s} F_n))$ in \eqref{e.prufer002}.
\begin{lemma}
For any $1\leq s\leq T$, $N\geq N_0(\varrho,\lambda)$,
\begin{equation}\label{e.LDT031}
\mathrm{mes}\left[(x,y)\in\mathbb{T}^2:\left|\frac{1}{N}\sum_{n=T+1}^N\mathrm{Re}(z^{2T-s}\zeta_{n-T}^2F_{n-s}F_n)\right|>\lambda^3 \right]\leq \mathrm{exp}(-C_\lambda N).
\end{equation}
\end{lemma}
\Proof
Let us consider two cases:\\

\textbf{Case 1.}  If $1\leq s\leq\frac{T}{2}$, set $\zeta_{n-T}^2=\xi_{n-T}+i\eta_{n-T}$. Then
it suffices to prove
\begin{equation}\label{e.LDT08}
\mathrm{mes}\left[(x,y)\in\mathbb{T}^2:\left|\frac{1}{N}\sum_{n=T+1}^{N}\xi_{n-T}(x,y)\mathrm{Re}(F_{n-s}F_n-\langle F_0F_s\rangle)\right|>\lambda^3\right]\leq\mathrm{exp}(-C_\lambda N).
\end{equation}
With $n=jT+r$, one has
\begin{equation*}
\begin{split}
&\frac{1}{N}\sum_{n=T+1}^{N}\mathrm{Re}(F_{n-s}F_n(x,y)-\langle F_0F_s\rangle)\xi_{n-T}(x,y)\\
&=\frac{1}{N}\sum_{r=1}^{s}\frac{1}{N/T}\sum_{j=1}^{[N/T]-1}\mathrm{Re}(F_{jT+r-s}F_{jT+r}-\langle F_0F_s\rangle)\xi_{(j-1)T+r}(x,y)+O(\frac{T}{N})\\
&=\frac{1}{N}\sum_{r=1}^{s}\frac{1}{N/T}\sum_{j=1}^{[N/T]-1}\mathrm{Re}(F_{jT+r-s}F_{jT+r}-\langle F_0F_s\rangle)[\underline{v}^{+} (\underline{v}^{+}\cdot (x,y)) ]\xi_{(j-1)T+r}(x,y)+O(\frac{T}{N}).
\end{split}
\end{equation*}
Let $f(\tilde{y})=\mathrm{Re}(F_{jT+r-s}F_{jT+r}-\langle F_0F_s\rangle)(\underline{v}^{+}\tilde{y})$, $g_{j,r}(\tilde{y},\tilde{y}')=\xi_{n-T}(\varrho^{T-n}\underline{v}^{+}\tilde{y}+\underline{v}^{-}\tilde{y}')$. Note that
\begin{equation*}
\begin{split}
&\mathrm{Re}(F_{jT+r-s}F_{jT+r}-\langle F_0F_s\rangle)(\varrho^{-r}\underline{v}^{+}\tilde{y})\xi_{(n-1)T+r}(\varrho^{-r}\underline{v}^{+}\tilde{y}+\underline{v}^{-}\tilde{y}')\\
&=f(L^j\tilde{y})g_{j,r}(L^{j-1}\underline{v}^{+}\tilde{y}+\underline{v}^{-}\tilde{y}')
\end{split}
\end{equation*}
with $L=\varrho^T \sim \lambda^{-100}$. Besides, $\left|\frac{\partial}{\partial\tilde{y}}g_{j,r}(\tilde{y},\tilde{y}')\right|\leq C\lambda$ and $hm_fh\leq \lambda^{-12}$. Hence, with $t\sim \lambda^3$,
\begin{equation*}
\begin{split}
{\rm LHS \ of\ } \eqref{e.LDT08}&\leq \exp(-t\lambda^3N/T)\int_{0}^1\sup_{|I|=1}\int_{I}\exp\left(t\sum_{j=1}^{[N/T]}f(L^j\tilde{y})g_{j,r}(L^{j-1}\tilde{y}+\tilde{y}')\right)d\tilde{y}d\tilde{y}'\\
&\leq \exp(-t\lambda^3N/T+Ct^2N/T)\leq \exp(-C_\lambda N).
\end{split}
\end{equation*}
Notice that condition \eqref{line02} is satisfied, since $(1+C\lambda)L^{-\frac{1}{2}}+\lambda^{-12}L^{-\frac{1}{2}}\leq t.$

\medskip
\textbf{Case 2.} If $\frac{T}{2}<s\leq T$, recall $\zeta_{n-T}^2=\xi_{n-T}+i\eta_{n-T}$. Set $f(\tilde{y})=\mathrm{Re}(F_n)(\underline{v}^{+}\tilde{y})$ and $g_{j,r}(\tilde{y},\tilde{y}')=\mathrm{Re}(F_{n-s}\xi_{n-T})(\varrho^{T-n}\underline{v}^{+}\tilde{y}+\underline{v}^{-}\tilde{y}')$, where $n=jT+r$. The estimation is the same as the above. Then we have
\begin{equation*}
\mathrm{Re}(F_{jT+r})(\varrho^{-r}\underline{v}^{+}\tilde{y})\mathrm{Re}(F_{jT+r}\xi_{(j-1)T+r})(\varrho^{-r}\underline{v}^{+}\tilde{y}+\underline{v}^{-}\tilde{y}')=f(L^j\tilde{y})g_{j,r}(L^{j-1}\tilde{y},\tilde{y}')
\end{equation*}
with $L=\varrho^{T-s+1} \sim \lambda^{-100}$. Moreover,
\begin{equation*}
\begin{split}
\left|\frac{\partial}{\partial\tilde{y}}g_{j,r}(\tilde{y},\tilde{y}')\right|&\leq (1+2\lambda)\varrho^{n-s}\varrho^{T-n}+C\lambda\varrho^{n-T-1}\varrho^{T-n}\\
&=(1+2\lambda)\varrho^{T-s}+C\lambda\varrho^{-1}\\
&\le (1+2\lambda)\varrho^{T-s+1}+C\lambda,
\end{split}
\end{equation*}
and up to a suitable truncation, $hm_f(h)\leq \lambda^{-12}$.
With $t\sim\lambda^2$, $[(1+2\lambda)\varrho^{T-s+1}+C\lambda]L^{-\frac{1}{2}}+\lambda^{-12}L^{-\frac{1}{2}} \leq t$,
 hence
\begin{equation*}
\begin{split}
\mathrm{mes}\left[(x,y)\in\mathbb{T}^2: \left|\frac{1}{N}\sum_{n=T+1}^{N}\mathrm{Re}(F_{n-s})\mathrm{Re}(F_n\xi_{n-T}(x,y))\right|>\lambda^3\right]<\exp(-C_\lambda N).
\end{split}
\end{equation*}
\qedbox

\medskip

Based on the above discussion, we can now start to prove Theorem \ref{mainth02}.

{\bf Proof of Theorem \ref{mainth02}} It will be enough to show that the large deviation estimate holds for $\frac{1}{2N}\log|\varphi_N(z)|^2$. We have
\begin{align}
\frac{1}{2N}&\log|\varphi_N|^2-\frac{\lambda^2}{2}\mathcal{J}(\eta)\notag\\
&=\frac{\lambda^2}{2N}\sum_{n=1}^{N}|F_n|^2-\langle \bar{F}_0F_0\rangle\label{eq001}\\
& -\frac{\lambda}{N} \sum_{n=T+1}^N \mathrm{Re}(z^T \zeta_{n-T} F_n )\label{eq002} \\
& +\frac{\lambda^2}{N}\frac{N-T}{N}\sum_{n=T+1}^N \sum_{s=1}^T \mathrm{Re}\{z^s(\bar{F}_{n-s} F_n-\langle \bar{F}_0F_s\rangle)\}\label{eq003} \\
& -\frac{\lambda^2}{N} \sum_{n=T+1}^N \sum_{s=1}^T \mathrm{Re}z^{2 T-s} (\zeta_{n-T}^2 F_{n-s} F_n)\label{eq004}\\
&-\frac{\lambda^2}{4N}\sum_{n=1}^{N}(\zeta_n^2F_n^2+\bar{\zeta}_n^2\bar{F}_n^2)+O\left(\lambda^3 \log ^2 \lambda^{-1}\right)\label{eq005}.
\end{align}
\medskip
Now we consider \eqref{eq003}. It is clear that
\begin{equation*}
\begin{split}
\lim_{N\rightarrow\infty}&\frac{\lambda^2}{N}\sum_{n=T}^{N}\sum_{s=1}^{T}z^s\bar{F}_{n-s}F_n=\lim_{N\rightarrow\infty}\frac{\lambda^2}{N}\sum_{n=T-s}^{N}\left(\sum_{s=1}^{T}z^s\bar{F}_{n}F_{n
+s}\right)\\
&=\lambda^2\int_{\mathbb{T}^2}\sum_{s=1}^{T}z^s\bar{F}_{0}F_{s}dxdy=\lambda^2\sum_{s=1}^{T}z^s\langle\bar{F}_0 F_s\rangle\\
&=\lambda^2\sum_{s=1}^{\infty}z^s\langle\bar{F}_0 F_s\rangle-\lambda^2 \sum_{s>T}z^s\int_{\mathbb{T}^2}\bar{F}_{0}F_{s}dxdy.
\end{split}
\end{equation*}
Let $f=|\bar{F_0}|$, $f_1=\mathrm{Re}\bar{F_0}$, $f_2=\mathrm{Im}\bar{F_0}$. Take $h=2^{-s}$. Then
\begin{equation*}
\begin{split}
&|\bar{F}_0-m_{f}(h)|\\
&\leq|\mathrm{Re}(\bar{F}_0)-m_{f_1}(h)|+|\mathrm{Im}(\bar{F}_0)-m_{f_1}(h)|\\
&=O(2^{-s})=O(\lambda^3),
\end{split}
\end{equation*}
and
\begin{equation*}
\lambda^2 \sum_{s>T}z^s\int_{\mathbb{T}^2}\bar{F}_{0}F_{s}dxdy=\lambda^2 \sum_{s>T} m_{f}(h)z^s\langle F_s\rangle+O(\lambda^3).
\end{equation*}
Due to $\langle F_s\rangle=0$, we have
\begin{equation*}
\lim_{N\rightarrow\infty}\frac{\lambda^2}{N}\sum_{n=T}^{N}\sum_{s=1}^{T}z^s\bar{F}_{n-s}F_n=\lambda^2\sum_{s=1}^{\infty}z^s\langle\bar{F}_0 F_s\rangle+O(\lambda^3).
\end{equation*}
Taking real parts,
\begin{equation}\label{e.real001}
\lim_{N\rightarrow\infty}\frac{\lambda^2}{N}\sum_{n=T}^{N}\sum_{s=1}^{T}\mathrm{Re}(z^s\bar{F}_{n-s}F_n)=\frac{\lambda^2}{2}\mathcal{J}(\eta)-\langle \bar{F}_0F_0\rangle+O(\lambda^3).
\end{equation}

One applies the large deviation we have made to the corresponding terms: \eqref{e.LDT02} to \eqref{eq001}; \eqref{e.LDT03} to \eqref{eq002}; \eqref{e.real001} to \eqref{eq003}; \eqref{e.LDT031}to \eqref{eq004}; \eqref{e.LDT04} to \eqref{eq005}. Then we have
\begin{equation}\label{e.LDT09}
\mathrm{mes}\left[(x,y)\in\mathbb{T}^2:\left|\frac{1}{2N}\log|\varphi_N|^2-\frac{\lambda^2}{2}\mathcal{J}(\eta)\right|>\lambda^3\right]<e^{-C_\lambda N}.
\end{equation}

Clearly, the exact same large deviation estimate applies for the $\psi_N(z)$, since the only change is $F_n\rightarrow-F_n$. For
 \begin{equation*}
L_N(z)= \frac{1}{2 N} \log \left(\left|\varphi_{N}(z)\right|^{2}+\left|\psi_{N}(z)\right|^{2}\right),
\end{equation*}
we have
\begin{equation*}
\begin{split}
\mathrm{mes}&[(x,y)\in\mathbb{T}^2:e^{L_N(z)2N}>e^{2N(\frac{\lambda^2}{2}\mathcal{J}(\eta)+ \lambda^3)}]\\
&=\mathrm{mes}[(x,y)\in\mathbb{T}^2:|\varphi_N(z)|^2+|\psi_N(z)|^2>e^{2N(\frac{\lambda^2}{2}\mathcal{J}(\eta)+ \lambda^3)}]\\
&\leq \mathrm{mes}[(x,y)\in\mathbb{T}^2:|\varphi_N(z)|^2>\frac{1}{2}e^{2N(\frac{\lambda^2}{2}\mathcal{J}(\eta)+ \lambda^3)}]\\
&\qquad+\mathrm{mes}[(x,y)\in\mathbb{T}^2:|\psi_N(z)|^2>\frac{1}{2}e^{2N(\frac{\lambda^2}{2}\mathcal{J}(\eta)+ \lambda^3)}],
\end{split}
\end{equation*}
which implies that \eqref{LDT01} holds. The large deviation estimates tell us that $\Omega_N$, the set of $(x,y)$ 's where $L_N(z)$ has a large deviation from $\frac{\lambda^2}{2}\mathcal{J}(\eta)$ is exponentially small. Altogether we have
\begin{equation*}
\begin{split}
L(z)&=\lim_{N\rightarrow\infty}L_{N}(z)\\
&=\lim_{N\rightarrow\infty}\int_{\mathbb{T}^2\setminus \Omega_N}\frac{1}{2N} \log \left(\left|\varphi_{N}(z)\right|^{2}+\left|\psi_{N}(z)\right|^{2}\right)dxdy\\
&=\frac{\lambda^2}{2}\mathcal{J}(\eta)+O(\lambda^3),
\end{split}
\end{equation*}
i.e., \eqref{lyapunov01}.
\qedbox

\section{Green's function estimates}

Define the unitary matrices
\begin{equation*}
\Theta_{n}=\left(
\begin{matrix}
\bar{\alpha}_n &\rho_n\\
\rho_n & -\alpha_n
\end{matrix}
\right).
\end{equation*}
Denote $\mathcal{L}_{+}$, $\mathcal{M}_{+}$ by
\begin{equation*}
\mathcal{L}_{+}=\left(
\begin{matrix}
\Theta_0 &~ & ~\\
~& \Theta_2 & ~\\
~ & ~& \ddots
\end{matrix}
\right),
\mathcal{M}_{+}=\left(
\begin{matrix}
\mathbf{1} &~ & ~\\
~& \Theta_1 & ~\\
~ & ~& \ddots
\end{matrix}
\right),
\end{equation*}
where $\mathbf{1}$ represents the $1\times1$ identity matrix.
The analogous factorization of $\mathcal{C}$ is given by $\mathcal{C}=\mathcal{L}_+\mathcal{M}_+$.

\medskip

Let $\mathcal{C}_{[a,b]}$ denote the restriction of an half-line CMV matrix to the finite interval $[a,b]$, $b>a>0$, defined by
 $$\mathcal{C}_{[a,b]}=P_{[a,b]}\mathcal{C}(P_{[a,b]})^{\ast},$$
where $P_{[a,b]}$ is the projection $\ell^{2}(\mathbb{Z}_+)\rightarrow \ell^{2}([a,b])$. $\mathcal{L}_{+,[a,b]}$ and $\mathcal{M} _{+,[a,b]}$ are defined similarly. More details can be found in \cite[Theorem 4.2.5]{B.Simon}.

{However, the matrix $\mathcal{C}_{[a,b]}$ will no longer be unitary due to the fact that $\alpha_a$ and $\alpha_b$ satisfy $|\alpha_a|<1$ and $|\alpha_b|<1$. Thus we state the modification of the boundary conditions briefly. For more details, see \cite[Section 3.2 and 3.3]{zhu01}.

With $\beta,\gamma \in \partial \mathbb{D}$, define the sequence of Verblunsky coefficients
\begin{equation*}
\tilde{\alpha}_n=
\begin{cases}
\beta, \qquad & n =a;\\
\gamma,& n =b;\\
\alpha_n,& n \notin  \{a,b\}.
\end{cases}
\end{equation*}

The corresponding operator is defined by $\tilde{\mathcal{C}}$ and we define
$\mathcal{C}^{\beta,\gamma}_{[a,b]}=P_{[a,b]}\tilde{\mathcal{C}}\left(P_{[a,b]}\right)^*$.
$\mathcal{C}^{\beta,\gamma}_{[a,b]}$ is unitary whenever $\beta,\gamma \in \partial\mathbb{D}$, which can be verified.
Then, for $z\in\mathbb{C}$, $\beta,\gamma\in\partial\mathbb{D}$, the polynomials are defined by
$$\Phi_{[a,b]}^{\beta,\gamma}(z):=\det{\left(z-\mathcal{C}_{[a,b]}^{\beta,\gamma}\right)}, \qquad
\varphi_{[a,b]}^{\beta,\gamma}(z):=(\rho_a\cdots\rho_b)^{-1}\Phi_{[a,b]}^{\beta,\gamma}(z).$$


By \cite[Section B.1]{zhu01}, one can get
\begin{equation}{\label{2.16}}
\left|G^{\beta,\gamma}_{\Lambda}(x,y;z)(n_1,n_2)\right|=\frac{1}{\rho_{n_2}}
\left|\frac{\varphi_{[0,n_{1}-1]}^{\beta,\cdot}(z)\varphi_{[n_2+1,N]}^{\cdot,\gamma}(z)}
{\varphi_{[0,N]}^{\beta,\gamma}(z)}\right|.
\end{equation}

In the Schr\"{o}dinger case, by restriction the eigenequation $(H-E)\xi=0$ to a finite interval $\Lambda=[a,b]$, one can get two boundary terms, which yields the identity $\xi(n)=-G_{\Lambda}^{E}(n,a)\xi(a-1)-G_{\Lambda}^{E}(n,b)\xi(b+1)$.

But in the CMV case, the analog of this formula depends on the parity of the endpoints of the finite interval. Concretely, if $\xi$ is a solution of difference equation $\mathcal{C}\xi=z\xi$,
we define
\begin{equation*}
\tilde{\xi}(a)=
\begin{cases}
(z\overline{\beta}-\alpha_a)\xi(a)-\rho_a\xi(a+1), \qquad & a\text{ is even,}\\
(z\alpha_a-\beta)\xi(a)+z\rho_a\xi(a+1),& a \text { is odd,}
\end{cases}
\end{equation*}
and
\begin{equation*}
\tilde{\xi}(b)=
\begin{cases}
(z\overline{\gamma}-\alpha_b)\xi(b)-\rho_b\xi(b-1), \qquad & b\text{ is even,}\\
(z\alpha_b-\gamma)\xi(b)+z\rho_{b-1}\xi(b-1),& b \text { is odd}.
\end{cases}
\end{equation*}
Then for $a<n<b$, we have
\begin{equation}{\label{C1}}
\xi(n)=G_{[a,b]}^{\beta,\gamma}(n,a;z)\tilde{\xi}(a)+G_{{[a,b]}}^{\beta,\gamma}(n,b;z)\tilde{\xi}(b).
\end{equation}

Similarly to \cite[Proposition 13.2]{Bourgain02}, we have
\begin{lemma}\label{le.green01}
There exists $\Omega_1 \subset \mathbb{T}^2$ such that $\operatorname{mes}\left(\mathbb{T}^2 \backslash \Omega_1\right)=0$ with the following property: For every $(x,y)\in \Omega_1$ there is $n_0=n_0(x,y, \lambda)$ such that for all
$$
n>n_0(x,y, \lambda),z\in I_0 \text { and }0 \leq s_0 \leq n^8,
$$
one has
\begin{equation}\label{e.green3}
\left|L_n(z)-n^{-4} \sum_{s=1}^{n^4} \frac{1}{n} \log \left\|M_n\left(A^{s+s_0} (x,y), z\right)\right\|\right|<\lambda^3 .
\end{equation}
In particular, for every $(x,y) \in \Omega_1$ and $N>n_0^4(x,y, \lambda)$,
\begin{equation}\label{e.green4}
\left|G_{[0,N]}^{\beta,\gamma}(x,y;z)\left(n_1, n_2\right)\right|<\frac{\exp \left[\left(N-\left|n_1-n_2\right|\right) L(z)+O\left(\lambda^{3} N\right)\right]}{\left|\operatorname{det}\left[z-\mathcal{C}_{[0,N]}^{\beta,\gamma}(x,y;z)\right]\right|}
\end{equation}
for all $n_1, n_2 \in[0, N]$.
\end{lemma}

\medskip

Let $(x,y)\in \Omega_1$. In the view of \eqref{e.green4}, for any $N>n_0^4(x,y,\lambda)$ and $z \in I_0$,
\begin{equation}\label{e.green5}
\begin{split}
& \left|G_{[0,N]}^{\beta,\gamma}(x,y;z)\left(n_1, n_2\right)\right| \\
& \quad \leq \exp \left(-\left|n_1-n_2\right| L(z)+\left[N L_N(z)-\log \left\|M_N(x,y;z)\right\|\right]+O\left(\lambda^{3} N\right)\right).
\end{split}
\end{equation}
Combining \eqref{e.green5} with the large deviation theorem yields
\begin{equation}\label{e.green6}
\left|G_{[0,N]}^{\beta,\gamma}(x,y;z)\left(n_1, n_2\right)\right|<\exp \left(-\left|n_1-n_2\right| L(z)+C \lambda^{3} N\right),
\end{equation}
for all $n_1, n_2 \in[0, N]$ provided $(x,y)$ lies in a set whose complement is of measure $<e^{-C_\lambda N}$. Indeed, this set, which depends on $z$, consists of those $(x,y) \in \mathbb{T}^2$ such that
$$
L_N(z)<\frac{1}{N} \log \left\|M_N(x,y;z)\right\|+C \lambda^3.
$$

\section{The proof of Theorem \ref{mainth01}}
Recall that
\begin{equation}\label{e.F01}
\mathcal{J}(\eta)=\sum_{n=-\infty}^{\infty}e^{in\eta}\langle \bar{F}_0F_n\rangle.
\end{equation}
Fix some small constants $\delta>0$ and suppose that $I_0\subset[e^{i\delta},e^{i(\pi-\delta)}]\bigcup[e^{i(\pi+\delta)},e^{i(2\pi-\delta)}]$ is some nonempty interval, such that for any $z=e^{i\eta}\in I_0$,
\begin{equation}\label{e.spectral01}
|\mathcal{J}(\eta)|>c>0,
\end{equation}
where $c$ is some positive constant. According to the asymptotic formula of the Lyapunov exponent \eqref{lyapunov01}, there exists a constant $c_0>0$, such that
\begin{equation}\label{e.lyapunov111}
L(z)>c_0.
\end{equation}

Similar to \cite[Proposition 13.3]{Bourgain02}, we have
\begin{lemma}\label{le.LDT04}
Fix $\beta,\gamma\in\partial\mathbb{D}$, $\delta>0$ small. Let constant $N$ be an arbitrary large integer and $\lambda \in (\lambda_0,1)$, $\lambda_0\in(0,1)$. Define $S_N=S_N(\delta,\lambda)$ to be the set of those $(x,y)\in\mathbb{T}^2$ for some $N_1\sim N^{12}$, $z\in I_0$ and $N_1<k<2N_1$, where $\bar{N}=[e^{(\log N)^2}]$. The following conditions hold:
\begin{equation}\label{e.double01}
\left\|(\mathcal{C}^{\beta,\gamma}_{[0,N_1]}(x,y;z)-z)^{-1} \right\|>e^{N^2},
\end{equation}
\begin{equation}\label{e.double02}
\frac{1}{N}\log\|M_{N}(A^k(x,y);z)\|<L_{N}(z)-C\lambda^3.
\end{equation}
Then
\begin{equation}\label{e.double03}
\mathrm{mes}(S_N)\leq e^{-C_\lambda N}.
\end{equation}
\end{lemma}
Let
$$S:=\limsup_{N\rightarrow \infty}S_N,$$
where $S_N$ be as in Lemma \ref{le.LDT04},thus $\mathrm{mes}(S)=0$.

It is for all $(x,y)\in\Omega:=\Omega_1\setminus S$, $\Omega_1$ as in Lemma \ref{le.green01}, that we shall prove localization.
Fix some $\delta$, $\lambda$ and such a choice of $(x,y)$. According to the Schnol's theorem for CMV Matrices \cite[theorem 3.4]{Damanik01}, for a.e. $z\in \mathrm{spec}(\mathcal{C})$ with respect to spectral measure there exists $\xi$ such that
\begin{equation}\label{e.local01}
(\mathcal{C}-z)\xi=0,
\end{equation}
\begin{equation}\label{e.local02}
|\xi_0|=1,
\end{equation}
\begin{equation}\label{e.local03}
|\xi_n|<n^C   ~\mathrm{for~ all}~ n=1,2,\cdots,
\end{equation}
where $C$ is some positive integer. It therefore suffices to show that for any $z$ any polynomially bounded solution decays exponentially.
To this end, fix some large $N$. By our choice of $(x,y)$, if there is $N_1\sim N^{12}$ such that
\begin{equation}\label{e.local04}
||(\mathcal{C}^{\beta,\gamma}_{[0,N_1]}(x,y;z)-z)^{-1}||>e^{N^2}
\end{equation}
holds for $\beta,\gamma\in\partial\mathbb{D}$, according to Lemma \ref{le.LDT04}, \eqref{e.double02} fails for every $N_1<k<2N_1$, i.e.,
\begin{equation}\label{e.localization03}
\frac{1}{N_1}\log\|M_{N_1}(A^k(x,y);z)\|>L_{N_1}(z)-C\lambda^3.
\end{equation}

As usual, let
\begin{equation*}
G^{\beta,\gamma}_{\Lambda}(x,y;z):=(\mathcal{C}^{\beta,\gamma}_{\Lambda}(x,y;z)-z)^{-1}
\end{equation*}
be the Green's function. Consider intervals
$$\Lambda=\left[k-\frac{N^3}{2},k+\frac{N^3}{2}\right]=[a,b],$$
where $N_1<k<2N_1$.
By definition of $G_{\Lambda}^{\beta, \gamma}$ and because of \eqref{e.local03}, it will suffice to prove that for fixed constant $c_{1}>0$,
\begin{equation}\label{e.local05}
\max_{n_1, n_2\in\Lambda}\left|G^{\beta,\gamma}_{\Lambda}(x,y;z)(n_1,n_2)\right|<e^{-c_1|n_1-n_2|+O(\lambda^3N^3)}.
\end{equation}
Once \eqref{e.local05} holds, then one can get the required exponential decay property as follows,
\begin{equation*}
\begin{split}
|\xi_{k}| <&|G^{\beta,\gamma}_{\Lambda}(k,a)||\xi_{a}|+|G^{\beta,\gamma}_{\Lambda}(k,b)||\xi_{b}|\\
\lesssim & CN^{12C}e^{-c_{1}\frac{N^3}{2}}\\
\lesssim & e^{-N^2},
\end{split}
\end{equation*}
for large $N$, which is the required exponential decay property. According to Lemma \ref{le.green01} and \eqref{e.localization03}, \eqref{e.local05} follows.

\medskip

According to the above statement, we first need to prove that for $N_{1}\sim N^{12}$,
$$\|G^{\beta,\gamma}_{[0,N_{1}]}(x,y;z)\|>e^{N^2}.$$
Recalling $\xi_0=1$, \cite[Lemma 3.9]{Kr} implies
\begin{equation*}
1\leq\|\xi_{[0,N_1]}\|\leq2\|G_{[0,N_1]}^{\beta,\gamma}(x,y;z)\|(|\xi_{N_1+1}|+|\xi_{N_1}|)
\end{equation*}
It therefore suffices to prove
\begin{equation}\label{e.local07}
|\xi_{N_1+1}|+|\xi_{N_1}|\lesssim e^{-N^2}.
\end{equation}

Now, we give the estimation of $|\xi_{N_{1}}|$. Suppose for some $k\sim N^{12}$, $G_{[k,N_1-1+k]}^{\beta,\gamma}$ satisfies the estimate
\begin{equation}\label{e.J01}
|G_{[k,N_1-1+k]}^{\beta,\gamma}(x,y;z)(n_1,n_2)|<e^{-L(z)|n_1-n_2|+O(N_1^3\lambda^3)},
\end{equation}
which implies
\begin{equation*}
\begin{split}
|\xi_{N_1}| <& N^{12C} |G^{\beta,\gamma}_{[k,N_1-1+k]}(N_1,k)|+N^{12C}|G^{\beta,\gamma}_{[k,N_1-1+k]}(N_1,N_1-1+k)|\\
\lesssim &N^{12C} e^{-L(z)\frac{N_1}{2}+O(N_1^3\lambda^3)} \\
\lesssim & e^{-N^2}.
\end{split}
\end{equation*}
Obviously, we also have $|\xi_{N_1+1}|\lesssim e^{-N^2}$.  Now, we show how to ensure By \eqref{e.J01},
\begin{equation*}
\begin{aligned}
&\left|G_{\Lambda}^{\beta,\gamma}(x,y;z)\left(n_1, n_2\right)\right|\\
&= \left|G^{\beta,\gamma}_{N^3}\left(A^{k-\frac{N^3}{2}}(x,y)\right)\left(n_1, n_2\right)\right| \\
< & \exp \left(-c \left|n_1-n_2\right|L(z)+O\left(\lambda^3 N^3\right)+N^3\left[L_{N^3}(z)-N^{-3} \log \|M_{N^3}(A^{k-\frac{N^3}{2}} (x,y);z)\|\right]\right) .
\end{aligned}
\end{equation*}
Thus, we should find $k\sim N^{12}$ for which
\begin{equation}\label{e.local09}
\left|L_{N^3}(z)-N^{-3} \log \|M_{N^3}(A^{k} (x,y);z)\|\right|<\lambda^3.
\end{equation}
Applying \eqref{e.green3} with $n=N^3$, $s_0=N^{12}$ gives
\begin{equation}\label{e.local10}
\left|L_{N^3}(z)-N^{-12} \sum_{k= N^{12}}^{2N^{12}}\frac{1}{N^3}\log \|M_{N^3}(A^{k} (x,y);z)\|\right|<\lambda^3.
\end{equation}
Thus \eqref{e.local09} follows from \eqref{e.local10} by averaging over $k$, which implies \eqref{e.J01} holds.
\qedbox
\vskip1cm

\section*{Acknowledgment}

Y.L. and D.P. were supported in part by NSFC (No.11571327, 11971059). S. G was supported in part by NSFC (No.12201588) and Shandong Provincial Natural Science Foundation (No.ZR2022QA035) and Fundamental Research Funds for the Central Universities (No.202264007).

\vskip1cm



\begin{thebibliography}{[aa]}

\bibitem{Alon01} N.\ Alon, P.\ Erd\"{o}s, J.\ Spencer, \textit{The probabilistic method}, New-York: John Wiley and Sons, 1992.
\bibitem{Bour08} J.\ Bourgain, \textit{Positive Lyapounov exponents for most energies}, GAFA, 37-66, LNM 1745, Springer 2000.

\bibitem{Brin} M. Brin, G. Stuck, \textit{Introduction to Dynamical Systems}, Cambridge University Press, 2010.

\bibitem{Bour04} J.\ Bourgain, \textit{Green's Function Estimates for Lattice Schr\"{o}dinger Operators and Applications}, Annals of Mathematics Studies, 2004.

\bibitem{Bour07} J.\ Bourgain, \textit{Anderson Localization for Quasi-Periodic Lattice Schr\"{o}dinger  Operators on $\mathbb{Z}^d$, $d$ Arbitrary}, Geom. Funct. Anal., 17.3, 682-706, 2007.

\bibitem{Bourgain01} J.\ Bourgain, M.\ Goldstein, \textit{On nonperturbative localization with quasi-periodic potential}, Annals of Math., \textbf{152}, 835-879, 2000.


\bibitem{Bourgain03}  J.\ Bourgain, M.\ Goldstein,  W.\ Schlag, \textit{Anderson localization for Schr\"{o}dinger operators on $\mathbb{Z}$ with potentials given by the skew-shift}, Comm. Math. Phys., \textbf{220}, 583-621, 2001.
\bibitem{BGS02}   J.\ Bourgain, M.\ Goldstein,  W.\ Schlag, \textit{ Anderson localization for Schr\"{o}dinger operators on $\mathbb{Z}^2$ with quasi-periodic potential}, Acta Math., 188, 41-86, 2002.
\bibitem{Bourgain02} J.\ Bourgain, W.\ Schlag, \textit{Anderson localization for Schr\"{o}dinger operators on $\mathbb{Z}$ with strongly mixing potentials},  Comm. Math. Phys.,\textbf{215}, 143-175, 2000.
\bibitem{CW21} C. Cedzich, A. H. Werner, \textit{Anderson localization for electric quantum walks and skew-shift CMV matrices}, Comm. Math. Phys., \textbf{387}(3), 1257-1279, 2021.

\bibitem{CL} R.\ Carmona, J.\ Lacroix, \textit{Spectral theory of random Schr\"{o}dinger operators}, Boston: Birkh\"{a}user, 1990.

\bibitem{ChuSinai89} V.A. Chulaevsky, Ya. G. Sinai, \textit{Anderson Localization for the I-D Discrete Schr\"{o}dinger Operator with Two-Frequency Potential}, Comm. Math. Phys., \textbf{125}, 91-112, 1989.

\bibitem{C-S95} V.A. Chulaevsky, T. Spencer, \textit{Positive Lyapounov exponents for a class of deterministic potentials}, Comm. Math. Phys., \textbf{168}, 455-466, 1995.

\bibitem{CFKS}  H.L.\ Cycon, R.G.\ Froese, W. \ Kirsch, B.\ Simon,\textit{ Schr\"{o}dinger operators}. Berlin-Heidelberg-NewYork: Springer, 1987.

\bibitem{Damanik01} D. Damanik, J. Fillman, M. Lukic, W. Yessen, \textit{Characterizations of uniform hyperbolicity and spectra of CMV matrices}, Discrete Contin. Dyn. Syst. Ser. S., \textbf{9}, 1009-1023, 2016.

\bibitem{Davis-Simon} E.B. Davis, B. Simon, \textit{Eigenvalue estimates for non-normal matrices and the zeros of random orthogonal polynomials on the unit circle}, J. Approx. Theory \textbf{141},189-213, 2006.

\bibitem{F-P} A. Figotin, L. Pastur, \textit{Spectra of random and almost periodic operators}, Grundlehren der mathematischen Wissenshaften 297, Springer 1992.

\bibitem{Goldstein01} M. Goldstein, W. Schlag, \textit{H\"{o}lder continuity of the integrated density of the sates for quasiperodic Schr\"{o}dinger equations and averages of shifts of subharmonic functions},  Annals of Math., \textbf{154(1)}, 155-203, 2001.

\bibitem{Jito99}  S.Ya. Jitomirskaya, \textit{Metal-insulator transition for the almost Mathieu operator}, Annals of  Math., \textbf{150}(3), 1159-1175, 1999.

\bibitem{Klein05} S.\ Klein, \textit{Anderson localization for the discrete one-dimensional quasi-periodic Schr\"{o}dinger operator with potential defined by a Gevrey-class function}, J. Funct. Anal., \textbf{218}, 255-292, 2005.

\bibitem{Kr} H. Kr\"{u}ger, \textit{Orthogonal polynomials on the unit circle with Verblunsky coefficients defined by the skew-shift},  Int. Math. Res. Not., \textbf{18}, 4135-4169, 2013.

\bibitem{LPG} Y. Lin, D. Piao, S. Guo, \textit{Anderson localization for the quasi-periodic CMV matrices with Verblunsky coefficients defined by the skew-shift}, J. Funct. Anal., \textbf{285}, 109975, 2023.

\bibitem{B.Simon} B.\ Simon,  \textit{Orthogonal Polynomials on the Unit Circle. Part 1. Classical Theory}, AMS Colloquium Publications, vol. 54, American Mathematical Society, Province, RI, 2005.

\bibitem{B.Simon-2} B.\ Simon, \textit{Orthogonal Polynomials on the Unit Circle. Part 2. Spectral Theory}, AMS Colloquium Publications, vol. 54, American Mathematical Society, Province, RI, 2005

\bibitem{Sinai87} Ya. G. Sinai, \textit{Anderson Localization for One-Dimensional Difference Schr\"{o}dinger Operator with Quasiperiodic Potential}, J. Stat. Phys., 46.5-6, 861-909, 1987.

\bibitem{T07} N. Timothy, \textit{Positive Lyapunov exponents for a class of ergodic orthogonal polynomials on the unit circle}, J. Math. Anal. Appl., \textbf{327}, 977-990, 2007.

\bibitem{Wang02} F.\ Wang, \textit{A formula related to CMV matrices and Szeg\H{o} cocycles}, J. Math. Anal. Appl., \textbf{464}, 304-316, 2018.

\bibitem{Wang01} F.\ Wang, D.\ Damanik, \textit{Anderson localization for quasi-periodic CMV matrices and quantum walks}, J. Funct. Anal., \textbf{276}, 1978-2006, 2019.

\bibitem{zhu01} X.\ Zhu, \textit{Localization for random CMV matrices}, J. Appr. Theory, \textbf{298}, 106008, 2024.

\end{thebibliography}
\end{document}